\documentclass[12pt,letterpaper]{article}
\usepackage{amsmath}
\usepackage{pdfpages}
\usepackage{graphicx}
\usepackage{adjustbox}
\usepackage{amsthm}
\usepackage{array}
\usepackage{relsize}
\usepackage{natbib}
\usepackage[colorlinks=true, linkcolor=blue, citecolor=blue]{hyperref}
\usepackage{caption}
\usepackage{subcaption}  
\usepackage{multirow}
\usepackage{tabularx}
\usepackage{mathtools}

\newcommand{\blind}{0}
\newcommand{\e}{\begin{equation}}
	\newcommand{\ee}{\end{equation}}
\newcommand{\interior}[1]{%
	{\kern0pt#1}^{\mathrm{o}}%
}
\def\b{\ensuremath\boldsymbol}
\def\mc{\ensuremath\mathcal}

\newtheorem{theorem}{Theorem}
\newtheorem{proposition}{Proposition}

\newtheorem{lemma}[theorem]{Lemma}
\newtheorem{corollary}{Corollary}
\newtheorem{definition}{Definition}

\newtheorem{remark}{Remark}

\renewcommand{\qedsymbol}{\rule{0.7em}{0.7em}}
\begin{document}
		\if0\blind
		{
			\title{\bf Robust Strongly Convergent M-Estimators Under Non-IID Assumption}
				\author{K. P.Chowdhury, Johns Hopkins University}
				
				
			\medskip
			\maketitle 
		} \fi
		\if1\blind
		{
			\bigskip
			\bigskip
			\bigskip
			\begin{center}
				\title{\bf Robust Strongly Convergent M-Estimator Under Non-IID Assumption}
				\runtitle{Robust Strongly Convergent M-Estimator Under Non-IID Assumption}
			\end{center}
			\medskip
		} \fi
		\begin{abstract}
			\noindent
			M-estimators for Generalized Linear Models are considered under minimal assumptions. Under these preliminaries, strong convergence of the estimators are discussed and an expansion of the estimating operators are given in the non-i.i.d. case with the i.i.d. case shown as a particular application. Various consequences of the results are discussed for binary and continuous models.  

\vspace{.2in}
			

		\end{abstract}
		%
\section{Introduction}
Consider any Generalized Linear Model (GLM), \e\b y_i = c(\b x_i)'\b \beta  + \epsilon_i, \ee 
where $ \epsilon_i $s follow some unknown distribution(s) and need not be independent of each other and the $ c(X_i) $ are some desired continuous functions of the observed design variables (covariates). We are interested in estimating the slope of such a model, $ \epsilon_i $ almost surely without the need for traditional i.i.d. assumptions or complex and time consuming learning algorithms, which are based on assumptions which may or may not hold for any particular iteration. When the complete data $ (\b y_i, X_i) $ are observed in the random sampling case, existing results ensure the slope estimates are asymptotically unbiased, consistent and Gaussian distributed. Where as in the doubly censored case there are a host of existing results \citep{zhang1996linear, ren1997regression, ren2003regression} in both the i.i.d. and non-i.i.d. cases. However, to the best of the author's knowledge no general non-i.i.d. method exists (whether censored or otherwise) when one or more of the assumptions above do not hold.

In this paper we consider GLMs when such i.i.d. assumptions need not hold, which makes asymptotic inferences impractical. We further ensure that we consider a construction which give results identical to existing constructs should traditional assumptions hold. Of particular interest is the convergence specifications under which scientific inferences are made for such models under varying assumptions on the underlying model, and how this affects estimation of strongly convergent unbiased parameters of interst for stochastic models.

As such, consider the latent variable and non-latent variable regression framework. Quite generally, in the existing latent variable formulations pointwise convergence is not guaranteed for the methods even under very strong parametric assumptions on the probability of success across such methodologies (\cite{CHOWDHURY2021101112}, \cite{SDSS2021}). Accordingly, let the probability of successes be $ F $ and $ F^* $ for binary regression and latent variable formulations respectively, which are necessarily unknown. As such, the forthcoming discussion will make clear that the link constraint holding for each observation is an absolutely crucial component needed for almost sure convergence to hold. This assertion is true no matter the estimation technique involved such as \cite{tanner1987calculation} or MLE since the characteristics of the underlying imposed topology are crucial for a properly defined linear operator on the relevant abelian group. Thus, in the forthcoming in Section \ref{mathfoundations}, I argue various propositions regarding the existing latent and nonlatent framework for baseline discrete outcomes along with some baseline definitions. In Section \ref{impossibility}, I discuss the impossibility of almost sure convergence in the existing framework in the absence of additional assumptions. In Section \ref{unified} I present the foundational mathematical results that makes possible the estimation of strongly convergent M-estimators in a non-i.i.d. setting through Theorem \ref{theo4}. In Section \ref{binomiallatent}, Section \ref{nonsigma} and Section \ref{existence} various other proofs are elaborated, which are needed for the main result.

\section{Mathematical Foundations for an Unifying Framework}\label{mathfoundations}

In the following I present the importance of the link condition holding for each observation and state some general results as to how it extends the current GLM framework very broadly. To be precise, I first prove under what circumstance there is equivalency of the current binary and latent variable formulations. Then to present a more unified framework I present some minimal topological definitions which will be needed for the remainder of the Section. I then give results which illuminate why the current GLM framework cannot give convergence results which are almost sure. In the impossibility theorem I show the necessary and sufficient conditions needed for almost sure convergence, and finally, I present an unified almost sure convergence methodology in a mathematically rigorous manner.  

To see that the current latent variable and the binary formulations need not give equivalent results, we need to consider multiple criteria. In particular, note that since asymmetric and symmetric DGPs induce different constraints on the latent probability of success, it can be used to give us our first result.

\begin{proposition}\label{prop0}
	Let $ F $ be a distribution for the Bernoulli probability of success and $ F^* $ the distribution for the latent probability of success. Then a necessary condition for equivalence of the Binary Regression and Latent Variable specifications is that $ F=F^* $. 
\end{proposition}

\begin{remark}
The proof of the proposition may be found in Section \ref{binomiallatent}, and highlights the importance of latent variables which are unique for comparison across widely used models in the statistical sciences.
\end{remark}

\begin{remark}The above result seems rather intuitive, since almost always the true distribution of success is unknown.
\end{remark}

\begin{remark}
	However, it is more difficult to see that indeed the convergence is not even guaranteed pointwise between a latent and nonlatent formulations even if the true distribution of the probability of success is somehow known and assumed to be the same for both the binary and latent variable formulations. 
\end{remark}

\begin{remark}
	To see this first note that the probability of success in the binary case is given by $ F(\b{x_i}'\b{\beta_i}) $ and it is assumed that the probability of failure is given by $ 1 -  F(\b{x_i}'\b{\beta_i}) $, which appears to be a reasonable conclusion. Yet the above results and the uniqueness of the Jordan Decomposition implies that this relationship need not hold even pointwise for the latent formulation!
\end{remark}

\begin{proposition}\label{prop6}
	The binary regression and latent variable formulations are equivalent if and only if $ c_1\nu^+ = h(F^*) = F $ and $ c_2\nu^- = 1 - F = (1 - h(F^*)) $ a.e. on the measureable space $ (\lambda, \Sigma) $ and $ h $ is a monotonic function of $ F^* $ with $ \{c_1, c_2\} \in \b{R}\setminus \{-\infty, \infty\} $.
\end{proposition}



\begin{remark}
	The proof of this result can again be found in Section \ref{binomiallatent}. The results have several far reaching consequences, and shows that latent and complete data methodologies need not be equivalent even under strong parametric assumptions such as $ F = F^* $. We must consider the probability of success and failure to be two separate measures for unique identifiability to hold. It further illuminates that no matter the estimation technique involved, simply assuming MLE results is not enough for congruence between the two models even in large samples, a finding readily validated in numerous empirical applications across the sciences (see for example \cite{CHOWDHURY2021101112} for a more detailed discussion).
\end{remark}

In fact, the result gives rise to several other relevant questions as to when the assumptions on the existing frameworks may and may not be supported. However, before discussing further result I highlight some topological definitions which will be needed going forward.

\subsection{Topological Definitions}

To facilitate this discussion the following definitions are asserted and may be found in almost any graduate level Topology book.

\begin{definition}
	A linear space $ \mc{X} $ is an abelian group with group operation addition, such that for a real number $ \alpha $ and $ \mathcal{u} \in \mc{X}\ and\ \{\alpha, \mathcal{\beta}\} \in \mathcal{R} $ a scalar product $ \alpha.\mathcal{u} \in \mc{X} $ and the following properties hold, \begin{itemize}
		\item $ (\alpha+\mathcal{\beta}).\mathcal{u} = \alpha.\mathcal{u} + \mathcal{\beta}.\mathcal{u} $.
		\item $ \alpha.(\mathcal{u}+\mathcal{v}) = \alpha.\mathcal{u} + \alpha.\mathcal{v};\ \mathcal{v} \in X $.
		\item $ (\alpha.\mathcal{\beta}).\mathcal{u} = \alpha.(\mathcal{\beta}.\mathcal{u})\ and\ 1.\mathcal{u} = \mathcal{u} $.
	\end{itemize}
\end{definition}

The addition and scalar multiplications are defined pointwise on all of $ \mc{X} $. On this space we may define a norm $ ||.|| $ as if $ \mathcal{u}, \mathcal{v} \in \mc{X} $ and $ \alpha \in \mathcal{R} $ then $ ||\mathcal{u}|| = 0,\ if\ and\ only\ if\ \mathcal{u} =0 $, $ ||\mathcal{u} + \mathcal{v}|| \le ||\mathcal{u}|| + ||\mathcal{v}|| $ and $||\alpha\mathcal{u}|| = |\alpha|||\mathcal{u}||$.

For two normed linear spaces we may define a linear operator as follows. \begin{definition}
	Let $\mathcal{X}$ and $ \mathcal{Y} $ be linear spaces. A mapping $ T: \mathcal{X} \rightarrow \mathcal{Y} $ is called a linear operator if for each $ \mathcal{u}, \mathcal{v} \in \mathcal{X} $ and real numbers $ \alpha $ and $ \mathcal{\beta} $,
	\e T(\alpha\mathcal{u}+\mathcal{\beta}\mathcal{v}) = \alpha T(\mathcal{u}) + \mathcal{\beta}T(\mathcal{v}). \ee
\end{definition}

In addition to the elementary definitions above I will work on a particular type of linear space the Banach spaces, and define it accordingly below.

\begin{definition}
	A Banach space is a complete normed linear space.
\end{definition}

To be complete we need one more topological concept, that of the Hausdorff Separation property.

\begin{definition}
	A topological space equipped with the Hausdorff Separation Property implies that any two points on the topological space can be separated by disjoint sets.
\end{definition}

These definitions allow us to highlight the results below.

\begin{proposition} \label{prop7}
	Let $ \mc X $ be a compact Hausdorff space, and $ \mc Y $ a compact subspace on which the link condition holds. Then under the assumptions of Proposition \ref{prop4}, the existing Latent Variable and Binary Regression frameworks are equivalent if and only if the underlying probability of success (failure) is symmetric around the origin.
\end{proposition}

\textbf{Proof.} \textbf{Case I:} Assume that the assumptions of Proposition \ref{prop4} holds. Then using Theorem \ref{theo0}, Proposition \ref{prop4} and Proposition \ref{prop5} we have that $ \nu^+ $ can be extended to all of $ \mc X $ through a linear functional $ \mc L $ on which \e \mc{L(X)} \le \nu^+(\mc X). \ee

Then from elementary functional analysis (\cite{plax2002}) we know that for $ C(\mc X) $ the space of continuous real-valued functions normed by the maximum norm every bounded linear functional $ \mc L $ on $ \mc X $ and $ f \in C(\mc X) $ we have that, \e \mc{L}(f) = \int_{\mc X} f \nu^+(dx), \ee

where $ \nu^+ $ belongs to $ C' $, the space of all finite signed measures. Consider the measure space $ (X, \Sigma, |\hat{\nu}|) $ as in Proposition \ref{prop4}. Then, \e |\mc L| = 1. \ee

But by Hahn-Banach we can extend this linear functional to all of $ \mc X $, and therefore, define \e \label{prbFail} \nu^- = F_0^* = 1 - \int_{\mc X} f \nu^+(dx),\ee

to get the desired assertion.

Conversely, now assume that \eqref{prbFail} holds. Then again by Hahn-Banach we have that there exists a linear funtional $ \mc L_1: \mc X \rightarrow \mc Y $ such that \e
|\mc{L}_1| = 1, 
\ee
and, 
\e \label{prbFail2} \nu^+ = F^* = 1 - \int_{\mc X} f \nu^-(dx).\ee

Thus, it remains to prove the necessity of symmetry around the origin.

But this specification makes clear the circumstances under which we can assume a symmetric distribution for the probability of success for the existing latent variable or binary regression framework, when, \e \label{prbFail3} \nu^- = F_0^* = 1 - \int_{\mc X} f \nu^+(dx) = \nu^+ = F^* = 1 - \int_{\mc X} f \nu^-(dx).\ee

Thus, we are done. 

\textbf{Case II:} Now assume that the conditions of Proposition \ref{prop3} holds. This scenario is considerably more convenient to deal with. Consider the case that \e |\nu^+| = |\nu^-|. \ee

Then again using Hahn-Banach we may extend a linear functional \e \mc{L}_2: \mc X \rightarrow \mc Y, \ee such that \eqref{prbFail3} holds. Thus, we are done. \qedsymbol

Therefore, the circumstances under which such an assumption is justified would under most circumstances be considered extremely restrictive and unlikely to represent the underlying stochastic process! Yet this mathematical formulation and line of reasoning provides even further striking results regarding the impossibility of almost sure convergence in any GLM framework for the existing specification. These results are summarized below.

\subsection{The Impossibility of Almost Sure Convergence in the Current Framework}\label{impossibility}

To show that the existing GLM framework does not guarantee almost sure convergence in general for any estimation technique it is necessary for us to consider linear operators between linear spaces relevant to both the underlying systematic component and for the link function. The following lemma puts this into more concrete terms.

\begin{lemma}\label{lemm1}
	Let $ \mathcal{X} $ be a finite dimensional linear space and consider $ \mathcal{Y} $ as the linear subspace on $ \eta = g(\mu) = \lambda $, as defined before as the link condition that ties the systematic component to the mean. Then there exists no unbounded linear operator T such that \e
	T:\mathcal{X} \rightarrow \mathcal{Y}
	\ee
	is continuous.
\end{lemma}

\textbf{Proof.} First note that, by construction $ \mc{X} $ is a finite dimensional linear space. Then for \newline $ \b\beta \in \mc{R}^n, n < \infty$, $ \lambda(\mc{X}, \b\beta) \in \mc{Y} \subset \mc{X}$ is a linear subspace of $ \mc{X} $. Further by construction of a GLM we know that $ T: \mc{X} \rightarrow \mc{Y} $ must exist.  Furthermore, by construction this linear operator must be unique for any given $ \b\beta $. Assume, T is unbounded. Then, for any $ \mc{u} \in \mc{X}$,
\e ||T(\mc{u})|| \geq M||\mc{u}||, \forall\ \mc{u} \in \mc{X}. \ee
By definition, $ \mc{u}_n \rightarrow \mc{u} $ then pointwise convergence implies $ \{T(\mc{u}_n)\} \rightarrow T(\mc{u}) $. Suppose T is continuous. WLOG consider $ \epsilon = 1 $ at $ \mc{u} = 0$. Then we should be able to pick a $ \delta $ such that $ ||T(\mc{u}) - T(0)||<1 $ with $ ||\mc{u}||<\delta $. But by assumption T is unbounded. Thus, there exists no $ M \ge 0$ such that \e ||T(\mc{u})|| < M||\mc{u}||, \forall\ \mc{u} \in \mc{X}. \ee
Therefore, no such $ \delta $ exists since no such M exists with $ \delta = (M||\mc{u}||)^{-1} $. Therefore, T is not continuous as needed.\qedsymbol 

The above results are instructive. It is not possible to find a continuous linear operator between the sample space and the link function if the link function can be either infinite dimensional with respect to the strong topology or it takes nonfinite values. Even if we assume the observed explanatory variables are a finite sample from a finite dimensional space, perhaps we may disregard the infinite dimensional case, but we cannot disregard any undefined values taken by such a linear operator, as this implies the operator is not continuous. If the operator is not continuous, many well known convergence results fail to hold regardless of whether Bayesian or Frequentist estimation methodology is used. Indeed this further implies that the results of \cite{tanner1987calculation} need not be continuous as there may exist a linear operator which need not be bounded. Consequently, in \cite{albert1993bayesian} we need not have unique convergence to the mean regardless of the MCMC method used to identify the posterior even if the observed data coincide with the strong assumptions mentioned previously.

Further note that in \cite{SDSS2021} I discussed that equivalency between the binomial and latent variable regression specification relies on the likelihood principle. That is, the two likelihoods must be proportional to each other for the almost sure convergence to hold between them. Accordingly, consider the Logistic link function, \e ln\left(\frac{p_i}{1-p_i}\right) = \lambda_i.\ee

As discussed in \cite{CHOWDHURY2021101112} the link function is indeed not bounded as \e \{p_i\} \rightarrow 1 \implies ln\left(\frac{p_i}{1-p_i}\right) \rightarrow \infty. \ee

Therefore, there can be no continuous linear operator between the sample space and the link field equipped with either the normal or hausdorff topology. The statement of the result above is actually rather more innocuous than perhaps its implications may initially indicate. Consider the following corollary as a direct implication of the results above. 

\begin{corollary}\label{coro3}
	Let $ \mc{X} $ and $ \mc{Y}$ be as in Lemma \ref{lemm1}. Then the Logistic and Probit formulations are equivalent in the sense of Birnbaum for any continuous or discrete GLM formulation. 
\end{corollary}

\textbf{Proof.} To see this rather surprising result first note the existence of a latent variable formulation is guaranteed by Proposition \ref{prop3} and Proposition \ref{prop4}. Further note that in Proposition \ref{prop5} I showed that a monotonic transformation of $ \nu^+ $ and $ \nu^- $ such that \begin{eqnarray} \nu^+ = h(F^*) = F\ and\\ \nu^- = 1 - h(F^*) = 1 -F, \end{eqnarray}
would result in the same inference using Birnbaum's Theorem. Consider a simple application using the density functions of the Logistic,\e
\b{f_1}(\b y,X,\b{\beta_1}, \sigma_1) = \frac{exp(-[\b y - X.\b\beta_1]/\sigma_1^2)}{\left(1+exp(-[\b y - X.\b\beta_1]/\sigma_1^2)\right)^2},
\ee
and that of the normal density, \e
\b{f_2}(\b y,X,\b \beta, \sigma) = \frac{1}{\sqrt{2\pi \sigma^2}} exp(-[\b y - X.\b \beta]/\sigma^2).
\ee

Thus, if $ \sigma = \sigma_1 $ the two exponential densities are proportional to each other, with an adequate monotonic transformation and without the imposition of a link constraint should give similar inference results in large samples under i.i.d. assumption even if the model fit and prediction results differ. Since the current Latent Variable framework imposes fixing the variance of the latent distribution for identification purposes, we may readily apply the constraint that $ \sigma_1 = \sigma $. But the Bayesian and Frequentist formulations are identical under prior restrictions using existing latent variable framework. Therefore, the statement of the corollary must hold in either formulation. Thus, the assertion of the corollary then readily follows. \qedsymbol

Indeed the result is validated across the sciences where (see for example \cite{albert1993bayesian}, \cite{cameron2010microeconometrics}) they state each model's parameters appear to be a constant multiple of the other. Specifically, $ \b{\beta}_{Logit} \approx 1.6\b{\beta}_{Probit} $ and seem to apply quite well in empirical applications. This then gives one of the more important results of this work.
\newline

\begin{corollary}
	Let $ \mc{X} $ and $ \mc{Y}$ be as in Lemma \ref{lemm1}. Then pointwise convergence is not guaranteed for any $ F^* = F $, and the linear functional \e T:\mc{X} \rightarrow \mc{Y}, \ee and the statement holds whether we use a Bayesian or Frequentist formulation.
\end{corollary}

\textbf{Proof.} Under the conditions of  Lemma \ref{lemm1} we have that no continuous linear functional exists between the sample space $ \mc{X} $ and $ \mc{Y} $, for any sequence $ \{c(X)\b{\beta}\} \rightarrow \lambda \in \mc{Y} $, if for some n, $ \nu^+(\lambda) = \infty $ or $ \nu^-(\lambda) = -\infty $. Therefore, by Proposition \ref{prop3} and Proposition \ref{prop4} using Corollary \ref{coro3} we have that no continuous functional exists such that pointwise $ T:\mc{X} \rightarrow \mc{Y} $ holds. 

Furthermore, by Proposition \ref{prop4} we know that a latent variable formulation in the discrete case can be extended to the continuous case. Therefore, the statement of the result holds for any GLM.

Observe, that this result is independent of the linear functional used and therefore it is independent of whether a Bayesian or Frequentist operator is used in the estimation process. \qedsymbol

The implications of this resuly, is more pernicious than may appear. The issue becomes apparent when the probability of success or failure is exactly equal to 1. This necessitates the $ \b{\beta} $'s to be numerically large for the link condition to hold explicitly, even if not considered in the estimation process. This in turn can result in non-convergence. This is often dealt with in practice, by throwing away a particular observation that may be causing estimation issues. While such an approach can lead to convergence to some parameter, it is not clear why such a parameter should be equivalent to the true likelihood under Birnbaum's theorem. In fact, the mathematical results state that they should not be. Another approach often taken is to start from multiple starting points to get the best model fit results. However, even this approach does not guarante convergence to the true parameter even pointwise, since there are infinitely many models and starting points that can be considered. Therefore, the need for a more rigorous functional analysis approach that is applicable across the sciences seems clear. Accordingly, in \cite{CHOWDHURY2021101112} I ensured that almost sure convergence can be achieved without facing the issues of an unbounded, discontinuous linear functional pointwise for some functional specification of $ \lambda(X, \b\beta) $. In \cite{SDSS2021}, I presented a methodology which ensured that if the latent and binary formulations are equivalent, then we may assure almost sure convergence of the parameters. In the current formulation, using some of the results of the previous works, and some of the insights presented above, I present a methodology which ensures convergence to the unique measure using Jordan Decomposition as in \cite{SDSS2021}, but in a far more general framework, for any continuous, bounded link specfication subject to the link constraint holding for each observation.
\section{An Unified Almost Sure Convergence Methodology}\label{unified}
The preceding sections laid the foundations for almost sure convergence to the true parameters of interest under the binary regression framework. They have done so by harnessing the advantages of both the latent variable and binomial regression case to overcome their respective disadvantages. Specifically, we know that the binary error can only take one of two values, either 0 or 1. On the other hand, for the latent variable formulation we know the error can be continuous. In the nonparametric section I set the mathematical foundations for the completely new robust methodology using the Jordan Decomposition Theorem for a signed measure. Those results showed it to be superior to existing methods with improvements to the parametric version under various settings. However, it still ensured for equivalency to the binomial regression framework that $ \nu^- = 1- F^* = 1 - F $. However, in the examples above I argued that this condition need not hold at all even if $ F^* = F $. Therefore, in this section I outline a methodology where we may relax this restrictive constraint. 

Accordingly, consider a likelihood function as follows,\e \label{uniflik}
L(X|\b{\beta}) = \b{c}(\b{\nu^+}(\b \lambda))^{(k)}(\b{\nu^-}(\b \lambda))^{(n-k)},\b c \in \mc{R},\ k = \sum_{i=1}^ny_i.
\ee
Note that such a likelihood function can be supported anytime we have independence between the two sample spaces over the cutoff point of 0. The framework of \cite{kass1989approximate} and \cite{SDSS2021} ensures that this can be done as an extension of the methodology presented in \cite{CHOWDHURY2021101112}. As in the nonparametric case, the current formulation also allows for this cutoff to be based on a normalized posterior probability such as the median. Using the ability to run the estimation algorithm over the proportional posterior then allows us to extend the nonparametric methodology in a more robust formulation.
Therefore the link condition holding pointwise takes the following formulation, \begin{eqnarray} 
	(\b{\nu^+})^{\alpha_1^*} = \b \lambda_{|S^+}(X, \b\beta),\\
	(\b{\nu^-})^{\alpha_2^*} = \b \lambda_{|S^-}(X, \b\beta).
\end{eqnarray}

This general framework then requires a substantially more intricate proof to guarantee almost sure convergence through Latent Adaptive Hierarchical EM Like (LAHEML) algorithm \citep{SDSS2021} extended to all measure spaces whether finite or $ \sigma- $finite. To see this the theorems below make this formulation clear in a mathematically rigorous way. In doing so, it adds to some well-known results in Real Analysis and Pure Mathematics.

\begin{theorem} \label{theo3}
	Let $ (X, \Sigma, \nu) $ be a finite measure space, with $ \nu $ finite a.e. Then under Latent Adaptive Hierarchical EM Like Formulation,\e \hat{\beta} \xrightarrow[\text{}]{\text{a.s.}} \beta \ee in $ L^p(X, \nu) $, where $ 1 \le p \le \infty $ and $ p = \infty $ represents the essentially bounded case.  
\end{theorem}

\textbf{Proof.} Consider Theorem \ref{theo1}, where almost sure convergence was asserted for $ \nu $ finite or finite a.e. on X. Thus, the case for $ 1 \le p < \infty $, is immediate. It remains then to show the case for $ p = \infty $. Thus, as in Theorem \ref{theo2}, let \e f_n(\lambda(\hat{\beta}^{(j)})|y^{(j)}, \alpha^{(j)^*}) := f_n^j. \ee Let $ \{f_n\} $ be the sequence of functions  on X for all j. Then $ \{f_n\} $ is bounded and finite by construction for each $ i \in \{1, 2,...,n_j\} $, where $ \cup_{i=1}E_i^j $ are the respective disjoint covering sets. Let $ \epsilon > 0, $ then there exists a $ \delta > 0 $ for $ E \in \Sigma $ such that, \e if\ \nu(E) < \delta\ then\ \int_E|f_n| < \epsilon,\ee follows straightforwardly from finiteness over $ E $. Therefore, by Dunford-Pettis Theorem (\cite{royden2010real}), $ f_n $ is weakly compact. Thus, by the Kantorovich Representation Theorem, there exists a linear functional, \e T_{\nu}: L^{\infty}(X, \nu) = \int_Xf\ d\nu \rightarrow \mathcal{R}, \ee where T is an isometric isomorphism of $ (X, \Sigma, \nu) $ on to the dual of $ L^{\infty}(X, \nu) $ (Ibid). Thus, we are done. \qed

\begin{remark}
	Using this result we may now give the most fundamental theorem of this work. That any functional specification may be represented by a linear functional through a finite or $ \sigma- $finite signed measure. As such it may be thought of as an extension of Riesz-Markov, and to the best of knowledge seems to be a new result in the literature.
\end{remark}

\begin{theorem}\label{theo4}
	Let X be a locally compact and Hausdorff topological space such that $ (X, \Sigma, \nu) $ is a $ \sigma- $finite measure space. Then under Latent Adaptive Hierarchical EM Like Formulation,\e \hat{\beta} \xrightarrow[\text{}]{\text{a.s.}} \beta \ee in $ L^p(X, \nu) $, where $ p \in \{1,...,\infty\}, $ where $ p = \infty $ represents the essentially bounded case. 
\end{theorem}

\textbf{Proof.} Consider the space of functions on X, $\mathcal{L}$ which are essentially bounded such that, there exists some $ M \ge 0$ with \e |f| \le M\ a.e.\ on\ X. \ee

Then there exists a linear functional such that $ 0 \le |f| \le 1. $ Following existing set up (\cite{royden2010real}) for $ L^P $ spaces, define $ L^P(X, \nu) $ to be the collection of $ [L] \in \mathbf{L} $, as the collection of extended real valued functions on X, which are finite a.e. on X. Thus, integrability implies measureability for all $ f \in \mc{L}. $ That $ L^P $ is a banach space is well extablished and it is stated without proof going forward.

Accordingly, consider the dual of X, $ X^*. $ By Alaoglu's Theorem the unit ball on $ X^* $ is weak-* compact. Let \e \bar{B}^*(1) = \{\psi_k \in X^* : |\psi_k - \psi| \le 1\}. \ee 

Fix $ \delta_k > 0 $, then by Alaoglu there are finitely many $ \psi_i's\ s.t.\ 0 \le i \le 1 $, \e |\psi_i - \psi_k| < \delta_k, \ee for some $ x \in X $ such that we may define sets of the form \e \label{XK}X_k = \{x, x_{n_k} \in X: |\psi_k(x_{n_k}) - \psi(x)| <\frac{1}{2^{n_k}}\}. \ee 

Further note that X is $ \sigma-finite$ such that there exists disjoint open sets with $ X =\cup_{k=1}^nE_k, $ where $ E_1 = \{\omega\} \cup X \sim \cup_{k=2}^n E_k$, with $\omega$ the one point Alexandroff compactification of $ E_1. $ Now each $ E_k $ is endowed with the subspace topology from X. Let $ \Lambda_k $ be a dense subset of $ (-1,1) $, and define on $ E_k $ a normally ascending collection of open subsets $  \interior{O}_{\lambda_k} $ of $ E_k $ with $ \lambda_k \in \Lambda_k. $ Let $ f_k: E_k \rightarrow \mathcal{R}, $ with $ f_k = 1 $ on $ E_k \sim \cup_{\lambda \in \Lambda} \interior{O}_{\lambda_k}$, and otherwise setting \e f_k(x) = inf\{\lambda_k \in \Lambda: x_n \in \interior{O}_{\lambda_k}\}. \ee Then \e f_k: E_k \rightarrow [-1, 1], \ee is continuous. 

Thus, for each $ E_k $ we may define a normally ascending collection of diadic rationals such that,
\e E_k \subseteq \interior{O}_{\lambda_1} \cup \interior{O}_{\lambda_2}...\interior{O}_{\lambda_{n_k}} \subseteq \bar{\interior{O}_{\lambda_1}} \cup \bar{\interior{O}_{\lambda_2}}...\bar{\interior{O}_{\lambda_{n_k}}}, \ee for some $ n_k \in \mathcal{N}.$ Consider the collection of functions on \e \label{lrestrict} \mathcal{L}_{|\interior{O}_{\lambda_1} \cup \interior{O}_{\lambda_2}...\interior{O}_{\lambda_{n_k}}}, \ee and consider the product topology on it. Then by the Tychenoff Product Theorem, it is also compact. 

Take a sequence $ \{x_{n_k}\} \in E_k $ such that $ \{x_{n_k}\} \rightarrow x \in E_k \subseteq X $. Then $ f_k(\{x_{n_k}\}) $ is discontinuous at most at countably many points. Thus, consider a sequence of functions $  f_{n_k}$ on \eqref{lrestrict} and diagonalize it such that for some $ n_k $, the neighborhood for \eqref{XK} for $ \frac{1}{2^{n_k}} $ is not continuous. Take $ n_{k+1} > n_k$. Continue this process untill it terminates to get an unique sequence of possibly disjoint normally ascending collection of open sets $ \interior{O^*}_{\lambda_k \in \Lambda_k}. $ If the construction is disjoint, we have identified a collection of disjoint open sets that cover $ E_k $. If not use the fact that X is Hausdorff and locally compact such that, \e
E_{n_k=j} = E_k \sim \cup_{n_k=1}^{j-1}E_{n_k} \implies E_k = \cup_{n_k=1}^j E_{n_k}, \ee is a disjoint collection of open sets. Since $ f $ is finite a.e. on X, it is also finite a.e. on $ E_k $, and therefore the diagonalization process is valid for each $ E_k $. 

Therefore, by construction the restriction of the measure space $ (X, \Sigma, \nu) $ to \e (X_{|\cap_{\lambda_k}\interior{O*}_{\lambda_k \in \Lambda_k}}, \Sigma_{|\cap_{\lambda_k}\interior{O*}_{\lambda_k \in \Lambda_k}}, \nu_{|\cap_{\lambda_k}\interior{O*}_{\lambda_k \in \Lambda_k}}), \ee is also a measure space, since by construction $ \nu(\interior{O}_1) < \infty $. Moreover it is finite and by compactness, a finitely additive measure space. Since by the continuity of measure, \e \nu(E_k)\le \sum_{\lambda_k \in \Lambda_k} \nu(\bar{\interior{O}}_{\lambda_k \in \Lambda_k}) < \infty. \ee

Thus, using Theorem \ref{theo3} we have that $ E_k $ may be covered by a countable collection of open sets $ \interior{O}_{\lambda_k \in \Lambda_k} $, whose closure $ \bar{\interior{O}}_{\lambda_k \in \Lambda_k} $ contains $ E_k, $ with the closure being compact and Hausdorff with respect to the weak-* topology on it, and the restricted measure space on it is bounded and finitely additive. Thus, we may define a linear operator, \e T_{\nu_{|\cup_{\lambda_k}\interior{O}_{\lambda_k \in \Lambda_k}}}(f_k) \rightarrow \int_{\cup_{\lambda_k}\interior{O}_{\lambda_k \in \Lambda_k}}f_k\ d\nu_{\interior{O}_{\lambda_k \in \Lambda_k}}\ for\ all\ f \in \mathcal{L}^{\infty}(E_k, \nu_{|\cup_{\lambda_k}\interior{O}_{\lambda_k \in \Lambda_k}}). \ee
	
	Then, T is an isometric isomorphism of the normed linear space \newline $ (X_{|\cup_{\lambda_k}\interior{O*}_{\lambda_k \in \Lambda_k}}, \Sigma_{|\cup_{\lambda_k}\interior{O*}_{\lambda_k \in \Lambda_k}}, \nu_{|\cup_{\lambda_k}\interior{O*}_{\lambda_k \in \Lambda_k}}), $ on to $ \mathcal{L}^{\infty}(E_k, \nu_{|\cup_{\lambda_k}\interior{O}_{\lambda_k \in \Lambda_k}}), $ by the Kantorovich Representation Theorem. Using Theorem \ref{theo2} and Dunford-Petis for each $ \delta_k > 0 $ we may define an $ \epsilon_k $ such that there exists an $ M_k > 0 $ with $ n_k > N_k \in \mathcal{N} $,
	\e
	\int_{\{x \in E_k: f_k(x) \ge M_k\}}|f_k| < \epsilon_k,\ for\ all\ n_k >  N_k \in \mathcal{N}.
	\ee
	
	But, \e
	\nu(E_k) \le \nu(\cup_{n_k=1}^{j(n_k)}\interior{O}_{\lambda_k \in \Lambda_k}) < \infty,
	\ee
	thus by Borel-Cantelli note that for each $ n_k $, by the countable monotonicity of $ \nu_{|\interior{O}_{\lambda_k \in \Lambda_k}}, $ \e
	\lim\limits_{n_k \rightarrow \infty}\nu\left(\cup_{\lambda_k=n_k}^\infty\interior{O}_{\lambda_k \in \Lambda_k}\right) < 1,
	\ee
	
	thus all but finitely many of the $ x's \in \cup_{\lambda_k=n_k}^\infty\interior{O}_{\lambda_k \in \Lambda_k}$ belong to finitely many of the $ \interior{O}_{\lambda_k} $'s. Then by Egoroff we may choose an $ n_k > N_k,$ such that there exists a collection of subsets $ E_{n_k} $ of $ E_k $ such that $ f_{n_k} \rightarrow f_k$ uniformly on $ E_{n_k} $ but $ \nu(E_k \sim E_{n_k}) < \epsilon_k. $
	
	Thus we may write,
	\begin{align*}
		\big\vert\int_X f d\nu - \int_{\cup_{k=1}^\infty E_k} f_k\big\vert < \big\vert\int_X f d\nu - \int_{\cup_{\lambda_k}\interior{O}_{\lambda_k \in \Lambda_k}} f_k \big\vert + \\ \big\vert\int_{\cup E_{n_k}} f_k -  \sum_{E_{n_k}} T_{\nu_{|\cup_{\lambda_k}\interior{O}_{\lambda_k \in \Lambda_k}}}(f_k)\big\vert + \big\vert\sum_{E_{n_k}} T_{\nu_{|\cup_{\lambda_k}\interior{O}_{\lambda_k \in \Lambda_k}}}(f_k) - \int_{\cup_{k=1}^\infty E_k} f_k\big\vert
	\end{align*}
	
	Taking the limit of $ E_k(\lambda_k) \rightarrow \infty $ as $ \lambda_k \rightarrow \infty $ gives us then,
	\e
	\big\vert\int_X f d\nu - \int_{\cup_{k=1}^\infty E_k} f_k\big\vert < 3\epsilon_k,
	\ee
	
	take $ \delta_k $ such that $ \epsilon_k < 1/3, $ answers the $ \epsilon $ challenge. Thus strong convergence immediately implies that we may define a sequence of random variables $ X_{n_k} $ from $ E_{n_k} $ onto $ [0,1) $ such that, \begin{align}
		\nu\left(\lim\limits_{n \rightarrow \infty}\cup_{n_k=1}^n\big\vert\nu^{-1}X_{n_k} - \nu^{-1}\int_{X_{|E_k}}f\ d\nu \big\vert\right)  \le \nonumber \\
		\nu\left(\lim\limits_{n \rightarrow \infty}\cup_{n_k=1}^n\nu^{-1}\big\vert X_{n_k} - \int_{X_{|E_k}}f\ d\nu \big\vert\right) \le \nonumber \\
		\nu\left(\lim\limits_{n \rightarrow \infty}\cup_{n_k=1}^n\nu^{-1}\big\vert X_{n_k}-(X_{|E_k})\big\vert\right)  \le \nonumber \\ \left(\lim\limits_{n \rightarrow \infty}\sum_{n_k=1}^n\nu_{|E_k} \big\vert E_{n_k} - E_{n_k}\cap \interior{O}_{\lambda_k \in \Lambda_k} \big\vert\right). \end{align}
	
	Therefore, take $ N_{\bar{K}} = max(N_k) $ over all k such that if $ n_k = n(\lambda_k) > N_{\bar{K}} $ we get \e \nu\left(\lim\limits_{n \rightarrow \infty}\cup_{n_k=1}^n\big\vert\nu_k^{-1}X_{n_k} - \nu^{-1}\int_{X_{|E_k}}f\ d\nu \big\vert\right) = \lim\limits_{n(\lambda_k) \rightarrow \infty} \sum_{n_k}\frac{1}{2^{n_k}} \rightarrow 1.\ee
	Thus, we are done. \qedsymbol
	
\begin{remark} The result has some important consequences for probabilistic models, since an application of an appropriate linear operator ensures almost sure convergence for all integrable functions over the sample space. 
	\end{remark}

\begin{remark} 
	Therefore, the parameter estimates in any GLM may be estimated almost surely as M-estimators over fintie or $ \sigma- $finite measure spaces. 
\end{remark}

\begin{remark}
	The formulation above is unqiue and the operator used ensures that convergence can occur uniquely without the need for independence or identical distribution assumption.
\end{remark}

\begin{remark}
	If the independent and identically distributed results hold then it is but one formulation out of the infinitely many possible formulations, and the results attained hold for this case as well.
\end{remark}

\begin{remark}
	No assumption is required here other than a Hausdorff topology on compact topological spaces. Any measure that maps a sample space to the reals is shown to have this property locally.
\end{remark}

Thus, this result shows that any functional form may be approximated locally through the use of an adequate linear operator. Such an operator is dicussed in \cite{CHOWDHURY2021101112}, \cite{SDSS2021} and \cite{Dukedissertation21}, under the name Latent Adaptive Hierarchical EM Like (LAHEML) algorithm, and the reader is kindly referred to those sources for further details. For the current purpose, the existence of such an operator in conjunctions with the theorems above shows that neither independence nor identical distribution assumption is needed for our chosen statistical methodology. As such the results present the only known methodology to the author that allows such flexibility under minimal assumptions of a Hausdorff topology over compact spaces.	
	
\section{Nonequivalency of Current Nonlatent and Latent Formulations}\label{binomiallatent}
To see that the current latent variable and the binary formulations need not give equivalent results, we need to consider multiple criteria. The first of these have already been alluded to before. In particular, note that since asymmetric and symmetric DGPs induce different constraints on the latent probability of success, it gives the result below.

\begin{proposition}\label{prop0}
	Let $ F $ be a distribution for the Bernoulli probability of success and $ F^* $ the distribution for the latent probability of success. Then a necessary condition for equivalence of the Binary Regression and Latent Variable specifications is that $ F=F^* $. 
\end{proposition}

\textbf{Proof.} First note that from before in the symmetric case, \e
p_i = F(\boldsymbol{x'_i}\boldsymbol{\beta_i}) = F^*[y_i^*>0] = F^*[-\epsilon_i < \boldsymbol{x'_i}\boldsymbol{\beta_i}] = F^*[\epsilon_i < \boldsymbol{x'_i}\boldsymbol{\beta_i}] = F^*[\boldsymbol{x'_i}\boldsymbol{\beta_i}]. 
\ee
Assume the assumption of the proposition does not hold. Then,\e
F^*[-\epsilon_i < \boldsymbol{x'_i}\boldsymbol{\beta_i}] \ne F^*[\epsilon < \boldsymbol{x'_i}\boldsymbol{\beta_i}], 
\ee
and we have by definition \e p_i = F(\b{x'_i}\b{\beta_i}). \ee
Let $ \b{\beta_i} $ be given, further we know by construction $ \b{x_i} $ are considered fixed. Then,
\e p_i = F(\b{x_i'}\b{\beta_i}) \ne F^*[\epsilon <\b{x_i}\b{\beta_i}], \ee
a contradiction to our hypothesis. The conclusion of the result then follows straightforwardly. \qedsymbol
\begin{proposition}\label{prop6}
	The binary regression and latent variable formulations are equivalent if and only if $ c_1\nu^+ = h(F^*) = F $ and $ c_2\nu^- = 1 - F = (1 - h(F^*)) $ a.e. on the measureable space $ (\lambda, \Sigma) $ and $ h $ is a monotonic function of $ F^* $ with $ \{c_1, c_2\} \in \b{R}\setminus \{-\infty, \infty\} $.
\end{proposition}

\textbf{Proof.} For the backward direction let, $ c_1 = c_2 = 1 $ and $ F = F^* $ but $ 1 - F \ne 1 - F^* $. Then the statement clearly does not hold. Since then if $ \nu^+ = 1- \nu^-$, we have the binary regression assumptions may hold for the Jordan Decomposition of the signed measures, yet the latent variable formulation does not equal it even pointwise. Thus, the backward negation is immediate.

For the forward direction, now assume $ F = h(F^*) $ and $ 1-F = 1- h(F^*) $ a.e. on the relevant measureable space. Then,  \begin{eqnarray}
	L(y_i|\b x_i, \b{\beta}_i) &= {F_i}(1-{F_i}),\\
	L(y^*_i|\b{x}_i, \b{\beta_i}) &= h({F^*_i}).\\
	\implies F \propto h(F^*)\ and\ 1 - F &= 1 - h(F^*)\\
	\implies L(y_i|\b{x}_i, \b{\beta_i}) \propto L(\b{y}^*_i|\b{x}_i, \b{\beta}_i)\ since\ h(F^*) \propto F\ pointwise.
\end{eqnarray}
Since the signed measure $ \nu $ can be decomposed into finite measures by proposition \ref{prop3} and \ref{prop4}, using Birnbaum's Theorem and the likelihood principle, we would arrive at the same inference for each methodology. The statement then is verified. \qedsymbol

The results are striking to this mathematician and gives several far reaching consequences. It states that the two methodologies need not be equivalent even under the assumption that $ F = F^* $, even without assuming any measure theoretic applications. We must consider the probability of success and failure to be two separate measures for unique identifiability to hold. It further illuminates that no matter the estimation technique involved, simply assuming MLE results is not enough for congruence between the two models even in large samples, a finding readily validated in numerous empirical applications across the sciences (see for example \cite{CHOWDHURY2021101112} for a more detailed discussion).

In fact, the result gives rise to several other relevant questions as to when the assumptions on the existing frameworks can and cannot be supported. The results below highlight these considerations.
	
\section{Estimating Non $ \sigma- $finite Measure Spaces}\label{nonsigma}

Some of the results presented below may be found in \cite{SDSS2021}, and are recreated here for completeness.

\begin{proposition}\label{prop3}
	Let $ (\lambda, \Sigma, \nu) $ be a measure space as above and $ \nu $ a finite signed measure on it. Then, \e |\bar{\nu}|(\Sigma) = \bar{\nu}^+(S^+)+ \bar{\nu}^-(S^-)\ee is a probability measure, where $ \{\bar{\nu}^+, \bar{\nu}^-\} $ are positive measures and $ S^+ $ is positive, but $ S^- $ is negative w.r.t. the signed measure $ \nu $.   
\end{proposition}
\textbf{Proof}. First note that since by construction $ \nu $ is finite, we must have that both $ \nu^+ $ and $ \nu^- $ must also be finite. Further $ \Sigma $ is a semiring, and thus using the Caratheodory-Hahn Theorem, we know  there exists a $ \mu:\Sigma \rightarrow [0,\infty)$. Then the Caratheodory measure defined as $ | \bar{\nu} | = \bar{\nu}^+(S^+)+ \bar{\nu}^-(S^-)  $ induced by $ \mu $ is an extension of $ \mu $ and $ | \bar{\nu} |$ is an unique extension and it is a finite measure. The statement of the proposition then is a straightforward consequence. \qedsymbol

\begin{proposition}\label{prop4}
	Let $ \left(\lambda, \Sigma\right) $ be a measureble space with $ \nu $ a measure which is neither finite or $ \sigma-finite $ such that $ \nu(\Sigma) = \infty\ or -\infty $ and let $ \left(\lambda, \Sigma, | \bar{\nu} |\right) $ be a $ \sigma-finite $ measure space. Then if the signed measure $ \nu $ takes one of the values of $ \{-\infty, \infty\} $ then either $ y=1 $ or $ y=0 $ for every observation w.r.t. the $ \sigma-finite $ measure.
\end{proposition}
\textbf{Proof.} By definition the signed measure does not take both values of $ -\infty $ or $ \infty $, and the labels of $ y=1 $ and $ y=0 $ are arbitrary. Thus, WLOG let $ \mu = -\infty $ over $ S^- $ where $ S^- $ is defined as before. Further, note that $ \nu^+ $ must be a finite measure and $ \nu^- $ a $ \sigma- $finite measure. Choose $ \epsilon $ and define, \e E^+_n = \{ y \in \Sigma|\nu^+_k(y) > (1-\epsilon)\phi^+(y)\ \forall\ k \ge n \}, \ee and its complement \e E^c_n = \{ y \in \Sigma|\nu^-_z(y) > (1-\epsilon)\phi^-(y)\ \forall\ z \ge t \}, \ee
where $ \phi^+ $ is a simple approximation of $ \nu^+ $ and $ \phi^- $ is a simple approximation to $ \nu^- $ w.r.t. the signed measure $ \nu $. Let us define $ M^+>0 $ then to be the maximum over all values of $ \{\phi^+, \phi^-\} $. Choose an index $ N = max\{N^+, N^-\} $ where $ N^+ $ is s.t. $ \nu^+(S^+ \sim E^+_n) < \epsilon$ for all $ n_1 \ge N^+  $ and $ N^- $ is chosen s.t. $ \nu^-(S^- \sim E^c_n) < \epsilon$ for all $ n_2 \ge N^-  $.
Then by additivity over domains of integrals, linearity and monotonicity and using Fatou's Lemma, we may define, \begin{eqnarray}
	| \bar{\nu} | \le lim\ inf \frac{\sum_{i=1}^n\nu^-(E^c_i)}{\int_{S^-}\phi_n^-\ d\nu} + lim\ inf \frac{\sum_{i=1}^n\nu^+(E^+_i)}{\int_{S^+}\phi_n^+\ d\nu}.\end{eqnarray}
Thus by Theorem \ref{theo0} we have,
\begin{eqnarray}
	| \bar{\nu} | \le 0+ lim\ inf \frac{\sum_{i=1}^n\nu^+(E^+_i)}{\int_{S^+}\phi_n^+\ d\nu}.
\end{eqnarray}
But $ | \nu | \in [0, \infty]$ by definition and $ \nu^+ \in [0, \infty) $ by construction. Thus,
\e | \bar{\nu} | \ge 0 + lim\ sup \frac{\sum_{i=1}^n\nu^+(E^+_i)}{\int_{S^+}\phi_n^+\ d\nu}. \ee
Therefore,
\e \label{prbDen} | \bar{\nu} | = \frac{\sum_{i=1}^n\nu^+(E^+_i)}{\int_{S^+}\phi_n^+\ d\nu} \le 1. \ee
Since $ \int_{S^+}\phi_n^+\ d\nu < \infty $ \eqref{prbDen} is well defined and the assertion follows. \qedsymbol \newline

\begin{proposition} \label{prop5}
	Consider a Hahn-Decomposition of the measure space $ (\lambda, \Sigma, \nu) $ into $ \{S^+, S^-\} $ as defined before, where the signed-measure $ \nu $ WLOG takes the value of $ -\infty $ but is not $ \sigma- $finite. Then there exists a linear functional $ \mathcal{L} $ which extends any measure $ \nu^+ $ over $ S^+ $ to all of $ L^p(Q, |\bar{\nu}|) $, with $ |\bar{\nu}| $ as in Proposition \ref{prop4}, and Q is the measureable space $  (\lambda, \Sigma)  $ with $ 1 \le p < \infty $.
\end{proposition}

\textbf{Proof.} This is a consequence of the Hahn-Banach Theorem. First note let f be a non-negative function on $ L^p $ restricted to $ S^+ $.  Let $ \psi $ be a Simple Function on $ \lambda $ the subspace of $ L^p(X_{|S^+}, \Sigma, |\bar{\nu}|) $. Let $ f $ be any bounded, continuous function on the subspace $ \lambda $ of X. Thus, from elementary measure theory we know by the simple approximation theorem that there exists a sequence $ \psi_n $ such that,\e
|\psi_n - f|^p \le 2^p.|f|^p\ on\ \lambda\ for\ all\ n.
\ee 
Further by construction $ |f|^p $ is integrable so by Lebesgue Dominated Convergence we know that $ \{\psi_n\} $ converges. Therefore the simple functions are dense and subadditive for the metric induced by the norm on $ L^p $. That $ \{\psi_n\} $ is positively homogeneous is straightforward and thus, the result is asserted without proof here. 

Therefore, by the Hahn-Banach Theorem we have that there exists a linear functional $ \mathcal{L} $ such that,
\e
\mathcal{L}(\lambda) \le \{\psi_n\}(\lambda),  
\ee

and further that it can be extended to all of X with the same norm. \qedsymbol

\begin{theorem}\label{theo0}
	Let $ (X, \Sigma, \mu) $ be any measure space with $ \{E_k\}^{\infty}_{k=1} \subseteq A \subset \Sigma $ a collection of measureable sets for which $ \mu(A) = \infty $. Then there exists some $ E_j $ where j is a countable collection of some k, such that \e \mu(\cup_jE_j) < \infty. \ee 
\end{theorem}

\textbf{Proof.} \textbf{Case I:} Let A consist of singleton sets of infinite measure. 

Then there is nothing to prove as \e \mu(\cup_{i=1}^{n} A_i \sim A_{j\ne i}) < \infty, \ee where $ A_j $ indicate the collection of singleton sets of infinite measure.

\textbf{Case II:} Let A be dense in $ \Sigma $. We seek to find a countable collection of measureable sets that have finite measure while satisfying the condition of the theorem.

By construction, \e \sum_{k=1}^{\infty}\mu(E_k) \le \infty. \ee 

Choose a number $ m_{ik} \in \b R\setminus\{\-\infty, \infty\} $, such that \e E_{ik} = {x\in A: \mu(E_{ik} < 1/m_{ik}) }. \ee
Then by the finite additivity and countable monotonicity of a measure, there exists a disjoint collection of sets, $ E_{ik} \subset E_k $, such that there exists an enumeration of $ E_{ik} $ to the natural numbers where the following condition holds,
\e \mu\left(  \cup_{n}\cap_{ik\ge n}E_{ik} \right) = \sum_{i=1}^n \mu(E_{ik}) \le 1. \ee

Further, we know by Borel-Cantelli there exists measureable sets $ E_{jk} \subset E_k$ such that 
\e \mu(\cap_{n}\cup_{jk\ge n}E_{jk}) = 1.\ee

Thus, there exists a disjoint countable collection of measureable sets in A such that $ \mu(\cup_{jk=1}^nE_{jk}) = 1 $. Therefore, there exists a set $ E_k $ with measure 1. Now by assumption $ \mu(A) = \infty $, therefore there exists some $ n_2 $ and $ n_3 $ belonging to the naturals such that,
\e 
\mu\left( (\cup_{k=1}^{n_2}E_k)\cup_{n_3\ge n_2}^{\infty}(A \sim (\cup_{k=1}^{n_2}E_k))\right) \le \mu(A).
\ee

Take $ n_2 \rightarrow\infty $ to get by finite additivity of measure that,
\e
\lim\limits_{n_2 \rightarrow \infty} \sum_{n_2}\mu(E_k) + \sum_{n_3\ge n_2}(A \sim (\cup_{k=1}^{n_2}E_k)) \le \infty,
\ee
thus, 
\e
\lim\limits_{n_2 \rightarrow \infty} \sum_{n_2}\mu(E_k) \rightarrow \infty. 
\ee

But k was arbitrary, and therefore denote by $ C = \cup_{k=1}^{j}E_k \subset A $ to get \e \lim\limits_{(k \rightarrow j)} \mu(C) \rightarrow \infty \implies \mu(C \sim A) < \infty,\ee

as needed.  \qedsymbol

\section{Existence and Uniqueness of Link Function Formulation} \label{existence}

Existence and uniqueness of the parameter estimates in statistical models then can be found below. While some of these results may be found in \cite{SDSS2021} and \cite{Dukedissertation21}, they are modified and applied here under the more general measure spaces discussed using Theorem \ref{theo4}.

\begin{theorem}\label{theo1}
	Let $ (\lambda, \Sigma) $ be a measureable space. Then there is an unique solution to any link modification problem, where the link constraint holds with equality in the Generalized Linear Model Framework for some $ \alpha^* \in R\setminus \{-\infty,\infty\} $, given $ \hat{F}_i $ $ \notin \{0,1\} $, $ X \notin \{0, \infty, -\infty\} $ element wise for each $ i  \in \{1,...,n\}$ and $\b \beta_j \notin \{\infty, -\infty\}$ with $ j \in \{1,...,(k+1)\} $.
\end{theorem}
\textbf{Proof.} \textit{Case I: $ (\lambda, \Sigma, | \nu |) $ is a finite measure space with $ | \nu | $ finitely additive and countably monotone.}

First note that by construction, 
\e
\{ | \nu |, \lambda(\mathbf{x}, \beta)\} < \infty \implies | \nu |^{\alpha^*} = \lambda(\mathbf{x}, \beta),  
\ee 

holds for some $ \alpha^* \in R\setminus \{-\infty, \infty\}$ by the density of the reals since,\e
|| \nu || = sup \sum_{i=1}^n | \nu |(E_i) < \infty.
\ee 

Thus, $ || \nu || $ is of bounded variation and $ | \nu | $ may be represented as the difference of two monotonic functions. As such, there exists a function $ g \in R\setminus \{-\infty, \infty\} $ such that  \e
g_{| \nu |} = | \nu |[I] = \hat{\textbf{F}}
\ee
where $ I $ is any countable collection of measureable sets covering $ R\setminus \{-\infty, \infty\} $. Such a covering exists from the compactness of the support on $ | \nu | $ and the assertion follows.

\textit{Case II: $ (\lambda, \Sigma, \nu) $ is a signed measure space.}

\hspace{0.2in}\textit{Subcase A: $ \nu $ is finite a.e. on $ \Sigma $.}

In this case we are back at \textit{Case I} above and the results hold.

\hspace{0.2in}\textit{Subcase B: $ \nu $ is not finite a.e. on $ \Sigma $.}

Consider the Caratheodory-Hahn extension to the measuareable space $ (\lambda, \Sigma, | \bar{\nu} |) $. From Proposition \ref{prop3} and Proposition \ref{prop4}, we know such an extension exists for $ \lambda $ lebesgue measureable. Consequently, using the results of \textit{Case I} again we arrive at the desired conclusion. \qedsymbol

This result ensures that an almost sure convergent methodology may be implemented under a particular linear operator.

\begin{theorem}\label{theo2}
	Given $ \alpha^* \in R\setminus \{-\infty,\infty\} $, and $ \hat{F}_i $ $ \notin \{0,1\} $, $ X  \notin \{0, \infty, -\infty\} $ elementwise for each $ i  \in \{1,...,n\}$ and $\b \beta_j \notin \{\infty, -\infty\}$ with $ j \in \{1,...,(J+1)\} $ subject to the link constraint holding for each observation, \begin{equation}
		\hat{\beta} \xrightarrow[]{a.s.} \beta.
	\end{equation}
\end{theorem}
\textbf{Proof.} Consider a Markovian framework and the folllowing cases below.

\textbf{Case I:} Let $ (\lambda, \Sigma, \nu) $ be a finite measure space with $ \nu $ a signed measure. 
\newline 
First note that the link condition holding for each observation implies following \cite{kass1989approximate} we may write for the appropriate measureable space $ (\lambda, \Sigma, | \nu |) $,\e p(y^*, \alpha^*, \beta|y) \propto p(y|\beta)f(y^*|\alpha^*, \beta, y)f(\alpha^*|\beta)f(\beta) \ee as the posterior distribution. Therefore, \e
p(\beta|y) \propto \int_{y^*}\int_{\alpha^*} p(y|\beta)f(y^*|\alpha^*, \beta, y)f(\alpha^*|\beta)f(\beta). 
\ee
Then\footnote{Following \cite{tanner1987calculation} I assume that $ \{y^*, \alpha^*\} $ both have compact support as the case for discrete support can be proved similarly.} considering all observations we may write,
\e p(\hat{\beta}^{(j)}|y)  \propto f_n\left(\lambda(\hat{\beta}^{(j)})|y^{(j)*}, \alpha^{(j)*}\right). \ee
For brevity I denote $ f_n\left(\lambda(\hat{\beta}^{(j)})|y^{(j)*}, \alpha^{(j)*}\right) $ as $ f^j_n $ going forward. Let T be an integral transform such that, \e
Tp(y|\beta)f(y^*|\alpha^*, \beta)f(\alpha^*|\beta)f(\beta) = p(\beta|y) \propto \int_{y^*}\int_{\alpha^*} p(y|\beta)f(y^*|\alpha^*, \beta)f(\alpha^*|\beta)f(\beta). 
\ee
Clearly, $ T $ is a bounded linear operator in $ L_1 $ by construction (the $ L_p $ case will be considered shortly) and define, \e f^{(j+1)}_n = (Tf^{(j)}_n)(\beta^{(j)}).
\ee
Then, define \e 
X_k = \{\lambda:  g^{(j)}_n =  
| f^{(j+1)}_n - (Tf^{(j)}_n)(\beta^{(j)}) | > k \}.
\ee
I claim $ g^{(j)}_n $ is an uniformly integrable sequence of functions w.r.t. $ | \nu | $ over $ \Sigma $. It is clear that by construction $ | \nu | $ is finite over $ \Sigma $. Thus, we may choose a natural number N such that if $ k \ge N $ we have,
\e\label{rel1} | \nu | (\cup_{k=N}^\infty E_k) < \frac{1}{N},
\ee
for each k. Let $ \tilde{N} $ be the maximum of these indicies such that for all k \eqref{rel1} holds. Further by the continuity of the measure we may choose a disjoint collection of such sets. Let $ g^{(j)}_{nk} $ be the restriction of $ g^{(j)}_n $ to $ E_k $. We know it is finite over $ E_k $, so by the simple approximation lemma there is a simple function $ g^{(j)}_{nk} $ such that, \e
\int_{E_k}g_n^{(j)} - \int_{E_k} g^{(j)}_{nk} < \epsilon/2.
\ee 
But $ g_n^{(j)} -  g^{(j)}_{nk} \ge 0 $ thus,
\e
\int_{\Sigma}g_n^{(j)} - \int_{\Sigma} g^{(j)}_{nk} \le \ \sum_{k=1}^{\infty} \int_{E_k} | g_n^{(j)} - \int_{E_k}g^{(j)}_{nk} | \le \tilde{N}|\nu|(\Sigma) + \epsilon/2.
\ee 
By letting $ | \nu |(\Sigma) = \epsilon/(2\tilde{N})  $ we get our desired result. Thus, $ g_n^{(j)} $ is uniformly integrable.
Further, by contruction the link condition holding for each observation implies $ g^{(j)}_n \xrightarrow[]{p.w.} | \nu |$ for each $ j $. Therefore, by the Vitali Convergence Theorem we have, 
\e 
\lim\limits_{j \rightarrow \infty} g^{(j)}_n = p(\beta|y).
\ee
%
%

Thus, we are done. \qedsymbol \newline



	\bibliographystyle{apalike}
	
	\bibliography{References.bib}
\end{document}